\documentclass[10pt]{article}
\usepackage{latexsym}
\usepackage{amssymb}
\usepackage{amsfonts}
\usepackage{amsmath}
\usepackage{theorem}
\usepackage{graphicx}
\usepackage{epstopdf}
\usepackage{mathrsfs}
\usepackage[breaklinks=true, linkbordercolor={0 0 1}]{hyperref}

\DeclareMathAlphabet{\mathpzc}{OT1}{pzc}{m}{it}

\usepackage{tikz}

\newtheorem{theorem}{Theorem}

\newtheorem{definition}{Definition}
\newtheorem{remark}{Remark}
\newtheorem{proposition}{Proposition}

\begin{document}

\title{Characterization of system signatures}

\author{Alessandro D'Andrea\footnote{Dipartimento di Matematica,  Sapienza Universit\`a di Roma,
P.le Aldo Moro 2, 00185 Roma, Italy;
{\tt<dandrea@mat.uniroma1.it>}},
\ Luca De Sanctis\footnote{Dipartimento di Matematica,  Sapienza Universit\`a di Roma,
P.le Aldo Moro 2, 00185 Roma, Italy;
{\tt<desanctis@mat.uniroma1.it>}}
}

\maketitle



\abstract
We suggest a purely combinatorial approach to a general problem in system
reliability. We show how to determine if a given vector can be the signature of
a system, and in the affirmative case exhibit such a system in terms on its structure function.
The method employs results from the theory of simplicial sets, and provides a full characterization
of signature vectors.


\section{Introduction}

The concept of signature of a system turns out to be very useful, as it provides
knowledge of the lifetime of a system in terms of its structure function
and single components' lifetime only. We refer to \cite{samaniego}
for a thorough introduction on the subject. In this paper, we give
a full characterization of signature vectors. This is performed by providing a criterion to check
whether a probability vector can be a signature. The method consists in numerical tests, 
and if the tests are positive it constructs explicitly 
the structure function a system with the given signature.
If a vector does not fulfill a certain technical requirement, the algorithm produces a family
which does not have the necessary algebraic properties so that no system can have 
the candidate vector as its signature.

The article is organized as follows.
In Section \ref{intro} we provide the main definitions and notions 
in the theory of system signatures.
In Section \ref{main} we recall the statement of the Kruskal-Katona theorem for
simplicial complexes and present the main result, i.e.
the algorithm that tests a vector to check whether it can be the signature of some system. 
The idea is based on the observation that a special family of sets of system components, 
which fully determines the system, has the same properties that are assumed in the
Kruskal-Katona theorem as to give place to a simplicial set. 
Then we simply import known theories and techniques from simplicial complexes into our context.
We conclude with some comments in Section \ref{conclusion}, where we also anticipate our next efforts,
part of them appearing in a forthcoming paper (\cite{dadss}).
Few introductory lines about simplicial complexes are listed in the final appendix,
to make the article self-contained.


\section{Basic notions}\label{intro}

In this section we recall some standard concepts and definitions from the theory of
system reliability (see \cite{samaniego} for more details).

Let $X=\{X_{1},\ldots , X_{n}\}$ be a set of $n\in\mathbb{N}$ continuous time binary
stochastic processes, interpreted as the state of the components of a system, 
each component $i$ being either down (or broken/off) or up (or working/on),
e.g. $X_{i}=0,1$ respectively, $i=1,\ldots , n$. We assume that all components are initially up
and when a component fails, it stays down forever.
A system deploys the components according to some design architecture and
is characterized by a \textit{structure function}
$\phi$ that indicates whether the whole system is up or down, 
for any given description of the states of individual components. 
In other words, the system may work even if some components are broken,
and given a subset $G\subseteq X$, interpreted as the set of working components,
the function
$$
\phi:2^{X}\to \{0,1\}
$$
tells us if the system is up ($\phi(G)=1$) or down ($\phi(G)=0$).
Common sense requires $\phi$ to be non-decreasing, which means
$A\subseteq B$ implies $\phi(A)\leq \phi(B)$. 
We assume that the random lifetimes of the components are continuous, in order to avoid ties in the
failures, and exchangeable. 
At the beginning all components (hence the whole system) work, and then
one at a time they fail (and stay broken), so that at some point the system stops working, say on
the $i$-th failure. The order in which components fail is a permutation $\sigma:\mathbb{N}\to
\mathbb{N}$
of the set $\{1,\ldots , n\}$, and this means that $\phi(\{X_{\sigma(i)},
\dots , X_{\sigma(n)}\})=1$ but $\phi(\{X_{\sigma(i+1)},\dots , X_{\sigma(n)}\})=0$.
We may rephrase this by saying that (for a given system $\phi$)
one and only one breakdown index $i\in\{1,\ldots, n\}$ is associated
with a given ordering of the failures (permutation) $\sigma$.
Let $N_{i}(\phi)$ be the number of permutations with breakdown index $i$, i.e.
such that $\phi(\{X_{\sigma(i)},\dots , X_{\sigma(n)}\})=1$
but $\phi(\{X_{\sigma(i+1)},\dots , X_{\sigma(n)}\})=0$.
Define $N(\phi)=(N_{1}(\phi),\dots, N_{n}(\phi))\in\mathbb{N}^{n}$.
\begin{definition}
The system \textit{signature} is the probability vector $s(\phi)=N(\phi)/n!$,
whose $i$-th entry $s_{i}(\phi)$ is the probability that the system stops working exactly
as the $i$-th failure of a component takes place.
\end{definition}
Let us introduce some standard terminology.
\begin{definition}
A subset $C\subseteq X$ is called a \textit{cut set} if the system is down when all the components
in $C$ are broken. A subset $P\subset X$ is called a \textit{path set} if the system is up whenever all
the components in $P$ work. 
A set of either type is said to be \textit{minimal} if none of its proper subsets
enjoy the same property.
\end{definition}
It is not difficult to see that the structure function is fully determined by the family of minimal
cut sets or equivalently by the family of minimal path sets. The system is thus completely defined
by its structure function or by its family of minimal cut or path sets. Minimality is the same as requiring that
there are no proper inclusions between elements of the family.
This one to one correspondence that maps non decreasing functions
$2^{X}\to \{0,1\}$ equal to zero at the empty set and equal to one at $X$ to
subsets of $2^{X}$ without proper inclusions
will be denoted by $\Omega$, and allows the following notation.
\begin{definition}\label{omega}
Given a structure function $\phi$, the corresponding
family of minimal cut sets is denoted by $\Omega(\phi)$. Similarly, for each family $\bar{\Omega}$
of subsets of $X$ without proper inclusions, there exists a unique structure function
$\phi_{\bar{\Omega}}$ whose family of minimal cut sets is $\bar{\Omega}$. 
\end{definition}
Clearly, $\phi_{\bar{\Omega}}=\Omega^{-1}(\bar{\Omega})$.

The next definition (\cite{samaniego}) is very convenient to capture the
relation between (minimal) cut and path sets in a system, by introducing some
sort of mirror system.
\begin{definition}
Given a system with structure function $\phi$, for all $A\subseteq X$ 
the \textit{dual} system has structure function $\phi^{*}$
$$
\phi^{*}(A)=1-\phi(X\setminus A)\ .
$$
\end{definition}
The map $\phi \mapsto \phi^{*}$ defines an involution in the set of structure functions on $X$.
Some well known remarks are in order. If $B^{\prime}\subset X$ is a minimal cut set and
$B\subset B^{\prime}$, then $X\setminus B$ is a path set;
conversely if $G^{\prime}\subset X$ is a minimal path set
and $G\subset G^{\prime}$, then $X\setminus G$ is a cut set.
Dual systems enjoy a remarkable property: a minimal
cut set for a system is a minimal path set for its dual system, and vice versa, since
duality is an involution. This correspondence
can be seen as a sort of time reversal, i.e. is yielded by reading the components of the
signature vector in reverse order. 
In fact, the signatures of two dual systems are related by
\begin{equation}
\label{duality}
s_{i}(\phi)=s_{n-i+1}(\phi^{*})
\end{equation}
for all $i=1,\ldots , n$.
Hence a vector can be a signature if and only if the same vector read in reverse order is also
a signature.


\section{Characterization of system signatures}\label{main}

In this section we present the main results of our study.

Let us start with a few preliminary observations.
Consider an ordering $\sigma$ of components' failures in the event that the
systems fails as the $i$-th failure of a component takes place.
This means that the set of the $i-1$ components that failed first
does not contain a cut set, and that the remaining $n-i+1$ components
include a path set (the system is still working at the time of the $i-1$-th failure)
and therefore a minimal path set as well.
We also know that the set of the $i$ components that failed first does contain
a cut set (and therefore a minimal cut set either), and the remaining $n-i$ components
do not include any path set. The component that gave place to the $i$-th failure
belongs to both a minimal cut set and a minimal path set. Since $\sigma$
indicate the order in which the components failed,
the component $X_{\sigma(i)}$ is the common element to the minimal cut and path
sets that appear in the first $i$ and last $n-i+1$ positions of the vector
$(X_{\sigma(1)},\ldots , X_{\sigma(n)})$ respectively.
This is a general fact.
\begin{remark}
Each minimal cut set intersects all minimal path sets, and the intersection consists of exactly 
one element. Conversely, each minimal path set intersects all minimal cut sets, and the 
intersection has cardinality one.
\end{remark}

The equivalence of structure functions and families of minimal cut or path sets
extends the duality introduced for structure functions to families of
minimal cut and path sets.
Given a family $\bar{\Omega}$ of subsets of $X$ without proper inclusions,
thanks to Definition \ref{omega},
we may define its dual family by
$$
\bar{\Omega}^{*}=\Omega(\phi^{*}_{\bar{\Omega}})
$$
regardless of whether it is interpreted as
family of cut or path set.
The family of minimal cut sets of a system
is also the family of the minimal path sets of the dual system (and vice versa, because duality is
an involution).
Duality is in essence the relation between minimal cut sets and minimal
path sets, which ultimately consists of a time reversal because of (\ref{duality}).
Translated in terms of structure functions this means
$$
\phi_{\bar{\Omega}^{*}}=\phi^{*}_{\bar{\Omega}}\ .
$$

We may use this identity to provide descriptions of dual families of minimal cut or path sets: 
let us do so in a fundamental example first.
If for some positive integer $l<n$ we choose $\Omega_{l}=\{\{1\},\{2\},\dots , \{l\}\}$,
then
$\Omega_{l}^{*}=\{\{1,2,\ldots , l\}\}$. This is the special case of series-parallel duality.
Now let $\phi_{l}=\phi_{\Omega_{l}}$. Using for simplicity the unnormalized signature
$N(\phi)=n!s(\phi)$, it is easy to figure that
$$
N_{i}(\phi_{l})=(n-l)!l! {n-i \choose l-1}\ , \ i<n-l+1
$$
and
$$
N_{i}(\phi^{*}_{l})=(n-l)!l! {i-1 \choose l-1}\ ,\ i>l\ .
$$
Let us list some explicit examples for $n=5$:
\begin{eqnarray*}
N(\phi_{1}) & = & (24,24,24,24,24) \\
N(\phi_{2}) & = & (48,36,24,12,0) \\
N(\phi_{3}) & = & (72,36,12,0,0) \\
N(\phi_{4}) & = & (96,24,0,0,0) \\
N(\phi_{5}) & = & (120,0,0,0,0).
\end{eqnarray*}
The corresponding duals are simply the same vectors with components in reverse order
because of (\ref{duality}).
The families of type $\Omega_{l}$ and $\Omega^{*}_{l}$ are very convenient: not only do they allow for simple calculations of their signatures, but they
make the signature of every system expressible as a $\mathbb{Z}$-linear affine
combination of their signatures. This is done using a standard inclusion-exclusion
procedure, and even though it may be difficult to perform explicit computations,
we obtain at least a representation theorem.
For example, if $\bar{\Omega}=\{\{1,2\},\{1,3\}\}$, then
$$
N(\phi_{\bar{\Omega}})=N(\phi_{\{\{1,2\}\}})+N(\phi_{\{\{1,3\}\}})-N(\phi_{\{\{1,2,3\}\}})\ .
$$
This is because the event ``the first $i$ failures involve components 1 and 2'' is the disjoint
union of the events ``the first $i$ failures involve components 1 and 2 but not 3'' and
``the first $i$ failures involve components 1 and 2 and also 3''. This last event is accounted
for even in the event ``the first $i$ failures involve components 1 and 3''. In other words,
$\{1,2,3\}$ is a superset of both $\{1,2\}$ and $\{1,3\}$ and is counted twice, so it has to be
subtracted once.

Since the system is exchangeable, $\phi_{\{\{1,2\}\}}$ and $\phi_{\{\{1,3\}\}}$
describe two equivalent systems, i.e. with same signature. 
Moreover, $\{\{1,2\}\} =\Omega^{*}_{2}$, $\{\{1,2,3\}\} =\Omega^{*}_{3}$,
and $\phi_{\Omega^{*}_{l}}=\phi^{*}_{\Omega_{l}}$. Therefore
\begin{eqnarray*}
N(\phi_{\bar{\Omega}}) & = & 2N(\phi^{*}_{2})-N(\phi^{*}_{3}) \\
{} & = & 2(0,12,24,36,48)-(0,0,12,36,72) \\
{} & = & (0,24,36,36,24) \\
s(\phi_{\bar{\Omega}}) & = & \frac{1}{5!}(0,24,36,36,24) =
\left(0,\frac{2}{10} , \frac{3}{10} , \frac{3}{10} , \frac{2}{10} \right)
\end{eqnarray*}
The general statement at this point follows quite naturally from the inclusion-exclusion principle, 
and it reads as follows.
\begin{theorem}
Let  $\bar{\Omega}$ be a family of subsets of $X$ without proper inclusions, interpreted
as the family of cut sets of an induced system with structure function $\phi_{\bar{\Omega}}$.
Then the signature of the system is
\begin{equation*}
s(\phi_{\bar{\Omega}})=\frac{1}{n!}
\sum_{\Gamma\subseteq\bar{\Omega}}(-1)^{|\Gamma|}s(\phi^{*}_{|\cup \Gamma|})
\end{equation*}
where $|\cdot|$ denotes the cardinality and $\cup\Gamma$ is the union of all the
sets in $\Gamma$.
\end{theorem}

We are now ready to describe an algorithm that produces a structure function inducing a given signature, thus providing a test to verify whether a vector arises as the
signature of a system, and a complete characterization of system
signatures. We start by introducing the Kruskal-Katona theorem.
Given two positive integers $n$ and $l$, it is known that there is a unique way to expand $n$ as a sum
of binomial coefficients as
\begin{equation*}
\label{expansion}
n=\binom{n_{l}}{l}+\binom{n_{l-1}}{l-1}+\cdots +\binom{n_{j}}{j}\ , 
\end{equation*}
with $n_{l}>n_{l-1}>\cdots > n_{j}\geq j\geq 1$.
Now for the given $n$ and $l$ define 
\begin{equation}
\label{plus}
n^{+}(l)=\binom{n_{l}}{l+1}+\binom{n_{l-1}}{l}+\cdots +\binom{n_{j}}{j+1}
\end{equation}
and
\begin{equation}
\label{minus}
n^{-}(l)=\binom{n_{l}}{l-1}+\binom{n_{l-1}}{l-2}+\cdots +\binom{n_{j}}{j-1}
\end{equation}
from the previous expansion. For the readers familiar with simplicial sets, 
these are the actions of the face and degeneracy maps. 
The next statement is a version of the theorem of Kruskal-Katona and 
offers a minimality constraint for simplicial complexes, 
with emphasis on the algebraic aspect of the sets composing the complex.
\begin{proposition}[Kruskal-Katona]\label{kk1}
Let $\underline{X}$ be a set of $n$ elements, $k$ and $l$ be given integers such that
$$
n\geq 1\ ,\ \ 1\leq l \leq n\ ,\ \ 1\leq k \leq \binom{n}{l}\ ,
$$
and let
$$
\mathcal{A}=\{A_{1},\ldots , A_{k}\}\ ,\ A_{i}\subseteq \underline{X}\ ,\ |A_{i}|=l\ ,\ i=1,\ldots , k\ .
$$
Let also
$$
\mathcal{A}^{-}=\{B:|B|=l-1\ ,\ \exists\ j : B\subset A_{j}\}\ .
$$
Then 
$$
\min_{\mathcal{A}}|\mathcal{A}^{-}|=n^{-}(l)\ ,
$$
where the minimum runs over all the families $\mathcal{A}$ of $k$ subsets of $\underline{X}$ of
cardinality $l$ and $n^{-}(l)$ is defined as in (\ref{minus}).
\end{proposition}
For the original proof and a more general analysis, see \cite{kruskal}.
The next statement is probably the most common version of the Kruskal-Katona theorem,
and provides a necessary and sufficient condition on the number of $l$-simplices for them to form
a complex. These numbers are the entries of the so called $f$-vector of the complex,
whose definition is recalled in Appendix \ref{app}.
\begin{proposition}[Kruskal-Katona]\label{kk2}
A vector $(f_{0},f_{1},\ldots , f_{d})$ is the $f$-vector of a simplicial $d$-complex if and only if
\begin{equation}\label{condition}
0\leq f^{-}_{l}(l)\leq f_{l-1}\ ,\ 1\leq l \leq d\ .
\end{equation}
Moreover, these two equivalent statements are equivalent to a third:
the subset $\Delta_{f}$ of the power set $2^{\mathbb{N}}$ consisting of the 
empty set together with the first (according to the reverse lexicographic order) $f_{l-1}$ subsets of 
$\mathbb{N}$ of cardinality $l$, for $l=1,\ldots , d$, is a simplicial complex. 
\end{proposition}
In our case $d=n-1$, and any total ordering can be chosen for the system components. 
Choosing initial segments (according to the reverse lexicographic
order) of size $f_{l}$ at level $l$ makes the number of elements al level $l-1$ minimal. 
There is a dual maximality condition which is equivalent to (\ref{condition})
$$
0\leq f_{l}\leq f^{+}_{l-1}(l)\ ,\ 1\leq l \leq d\ ,
$$
where $f^{+}_{l-1}(l)$ is defined as in (\ref{plus}).
For the original proofs and a more general analysis, 
see \cite{kruskal, katona}. The reverse lexicographic order simply 
reads backwards the strings, then sorts lexicographically. The advantage of considering the
reverse lexicographic order is that the list of the first (according to this order) $r\in\mathbb{N}$ elements
does not depend on the size of the alphabet (the size $n$ of the system, in our case).

Let us see how this applies in the context of system signatures. 
The algorithm that we are about to present is the translation of the proof of the Kruskal-Katona 
theorem where the role of $f$-vectors is played by the ``complement'' of the cumulative 
signature times the number of permutations of components, roughly speaking.
\begin{theorem}\label{sign}
Let $\bar{s}\in\mathbb{Q}^{n}$ be the candidate signature. Assume it is a probability vector.
For $l=1,\ldots , n$, define $f_{l}=\binom{n}{l}(\bar{s}_{l+1} + \cdots + \bar{s}_{n})$.

Then $\bar{s}$ is the signature of a system if and only if 
\begin{equation}
\label{test}
0\leq f^{-}_{l+1}(l)\leq f_{l}\ ,\ 1\leq l \leq n-1\ ,
\end{equation}
where $f^{-}_{l}(l)$ is defined as in (\ref{minus}).
\end{theorem}
\textbf{Proof.}
Let $\tilde{\Omega}$ be the family of cut sets of a system with structure function $\phi$, 
and consider $\mathcal{C}\equiv 2^{X}\setminus \tilde{\Omega}$.
Then $\mathcal{C}=\phi^{*-1}(0)$ is a simplicial complex. Similarly, using path sets instead of cut sets, 
$\mathcal{P}\equiv \phi^{-1}(0)$ is also a simplicial complex. 
(We will focus on cut sets only, but nothing changes if we consider path sets).
In fact, it is obvious that a superset of a cut or path set is still a cut or path set respectively.
Denote by $\mathcal{C}_{l}$ and $\tilde{\Omega}_{l}$ the set of elements of $\mathcal{C}$ 
and $\tilde{\Omega}$ respectively of cardinality $l$, so that $\mathcal{C}=\cup_{l}\mathcal{C}_{l}$
and $\tilde{\Omega}=\cup_{l}\tilde{\Omega}_{l}$. 
Clearly $\mathcal{C}_{l}\cup \tilde{\Omega}_{l}=\binom{X}{l}$ and 
$|\mathcal{C}_{l}|+|\tilde{\Omega}_{l}|=\binom{n}{l}$.
If $A\in \mathcal{C}_{l}$
and $x\in A$, then $A\setminus x \in \mathcal{C}_{l-1}$. Therefore the vector 
$(|\mathcal{C}_{1}|,\ldots , |\mathcal{C}_{n}|)$ is the $f$-vector of the simplicial complex $\mathcal{C}$.
Knowing this vector means knowing the family of cut sets too, since 
$\tilde{\Omega}_{l}=\binom{X}{l}\setminus \mathcal{C}_{l}$. But the vector 
$(|\tilde{\Omega}_{1}|,\ldots , |\tilde{\Omega}_{n}|)$ is the non-normalized cumulative signature
whose $l$-th component is $\binom{n}{l}S_{l}(\phi)=\binom{n}{l}(s_{1}(\phi)+\cdots + s_{l}(\phi))$.  
So given the candidate signature $\bar{s}$ we also know the non-normalized 
candidate cumulative signature $\bar{S}$
and hence the corresponding candidate $f$-vector for $\mathcal{C}$, 
with components $f_{l}=\binom{n}{l}(1-\bar{S}_{l})$, $l=1,\ldots , n$. Clearly
$1 - \bar{S}_{l} = \bar{s}_{l+1} + \cdots + \bar{s}_{n}$ by definition. 
The Kruskal-Katona theorem at this point provides a test to check whether 
such a vector is actually an $f$-vector, the test consists precisely of condition (\ref{test}),
as explained in Proposition \ref{kk2}.
$\Box$

Sometimes an equivalent procedure might be handier, especially for small systems.
Here it follows.
Let $\bar{N}\equiv n!\bar{s}\in\mathbb{N}^{n}$. 
\begin{enumerate}
  \item For each $l=1,\ldots , n$, sort in lexicographic order the subsets of $\{1,\ldots , n\}$
of cardinality $l$.
  \item Take the family $\Omega^{l}$ of the first
$(\bar{N}_{1}+\cdots + \bar{N}_{l})/(n-l)!l!=(\bar{s}_{1}+\cdots + \bar{s}_{l})\binom{n}{l}$ subsets,
with respect to the lexicographic order.
  \item Take the union $\cup_{l=1}^{n} \Omega^{l}$ of all the $\Omega^{l}$, $l=1,\ldots , n$,
  and extract the minimal family $\bar{\Omega}$.
  \item The function $\phi_{\bar{\Omega}}$ is the structure function of a system $X$ with $n$ 
  components and signature $\bar{s}$.
\end{enumerate}
Now considering Proposition \ref{kk1}, the same proof of theorem \ref{sign} also proves the following
test.

\textbf{Criterion.}
The family $\Omega^{l+1}$ should contain all the supersets (of cardinality $l+1$) of
at least one element from $\Omega^{l}$. 
If this is not the case, then the vector $\bar{s}$
cannot be the signature of a system since $\cup_{h}\Omega^{l}$ is not a simplicial complex.

This criterion is equivalent to Theorem \ref{sign}, and
the algorithm we presented is simply the Kruskal-Katona algorithm adjusted to 
work directly with the candidate non-normalized cumulative signature as opposed to 
its ``complementary'' vector with components $\binom{n}{l}(1-\bar{S}_{l}(\phi))$, $l=1,\ldots , n$.  
This is the reason why we sort strings lexicographically, because the
collection $\mathcal{C}$ is a simplicial complex, as opposed to $\tilde{\Omega}$. 
So instead of taking, as in the original Kruskal-Katona algorithm,
initial segments in each $\mathcal{C}_{l}$ according to reverse lexicographic order, 
we take final segments, i.e. initial segments according to the reverse ordering, 
which is the lexicographic order, in each $\tilde{\Omega}_{l}$.

We want to show that this second algorithm can be fairly fast in an explicit detailed example. 
In \cite{dadss}, we study two systems that are described by this vector.

\textbf{Example.}
Consider the vector $(0,3/10,2/5,3/10,0)$.
We pass easily to the non-normalized one $(0,36,48,36,0)$ by multiplying by $5!$. 

Start with $l=1$. We must take the first zero singletons.

Take $l=2$. We must take the first 36/12=3 subsets with two elements. These
are $\{1,2\}$, $\{1,3\}$, $\{1,4\}$.

Take $l=3$. We must take the first 84/12 = 7 subsets with three elements.
These are $\{1,2,3\}$, $\{1,2,4\}$, $\{1,2,5\}$, $\{1,3,4\}$, $\{1,3,5\}$, $\{1,4,5\}$, $\{2,3,4\}$.

Take $l=4$. We must take the first 120/84 = 5 subsets with four elements.
These are $\{1,2,3,4\}$, $\{1,2,3,5\}$, $\{1,2,4,5\}$, $\{1,3,4,5\}$, $\{2,3,4,5\}$.

Take $l=5$. We must take the first 120/120=1 subsets with five elements. This is
$\{1,2,3,4,5\}$.

From all these subsets we must extract a minimal family. It is not difficult to obtain
$\bar{\Omega}=\{\{1,2\},\{1,3\},\{1,4\},\{2,3,4\}\}$.

This fully determines the system,
and we can use the definition of minimal cut sets to determine the structure function
$\phi_{\bar{\Omega}}$ and verify that $N(\phi_{\bar{\Omega}})=(0,36,48,36,0)$.


\section{Conclusions and outlook}\label{conclusion}
We have introduced an algorithm, borrowed from the theory of simplicial complexes,
in the field of system reliability that checks if a given probability
vector can be a system signature, and in case constructs a system with that signature.
This completely characterizes the set of possible system signatures.
In a second paper (\cite{dadss}), we will show further results that follow from the analogy between
system signatures and $f$-vectors of simplicial complexes. Namely, we will show
that the only signature with first and last component both different from zero is the uniform one;
we will show that two systems with the same signature can be different (even up to permutation of
components). We will also show that a signature vector cannot have an isolated zero component,
and study the unimodal property of signatures. 

The bridge between probability in the theory of reliability and other fields where
the Kruskal-Katona theorem has proven to be fruitful, opens some promising
perspectives, as hopefully more than just the results of the current article and of \cite{dadss}.
We are investigating the possibility of a quantum theory of the signature, employing a
$q$-deformed binomial representation of integers, $q$-simplicial categories, etc..
We are also exploring ways to extend the study to include system availability,
as a second crucial quantity to evaluate in the context of RAMS (Reliability, Availability, Maintainability,
Safety) problems.


\appendix

\section{Simplicial complexes}\label{app}
In this section we recall some notions in the theory of simplicial complexes.

A \textit{simplicial complex} $K$ is a set of simplices such that
any face of a simplex from $K$ is also in $K$ and that
the intersection of any two simplices $\Sigma_{1},\Sigma_{2}\in K$ 
is a face of both  $\Sigma_{1}$ and $\Sigma_{2}$. A simplicial $d$-complex is a simplicial complex
where the largest dimension of any simplex in it is $d$. The $f$-vector of a simplicial $d$-complex
is the vector $(f_{0},f_{1},\ldots, f_{d})$ whose $i$-th component is 
the number of $(i-1)$-dimensional faces in the simplicial complex, and by convention $f_{0}=1$
unless the complex is empty.
The Kruskal-Katona theorem provides a full characterization of $f$-vectors
of simplicial complexes.

\section*{Acknowledgments}

The authors would like to thank Mario Marietti for making us aware that a former characterization
we had produced was actually equivalent to the one provided by Kruskal and Katona. 
The authors would like to express their gratitude to Fabio Spizzichino for suggesting the interesting
open problem they tackled here.
ADA was partially supported by ``La Sapienza'' Ateneo fundings.
LDS was partially supported by research funds of C. De Concini.


\end{document}